\documentclass{article}
\usepackage{amsmath,amssymb}

\begin{document}
\title{\bf Fuzzy $h$-ideals of hemirings}
 \normalsize

\author{{   Jianming Zhan}$^{a,*}$  , {Wies{\l}aw A. Dudek}$^b$  \\ {\small
 $^a$  {\em Department of Mathematics,
  Hubei Institute for Nationalities,}}\\ {\small {\it Enshi, Hubei Province,
   445000, P. R. China}}\\
  {\small  $^b$ {\em Institute of Mathematics and Computer Science,
         Wroc{\l}aw University of Technology,}}\\
         {\small {\em Wybrze\.ze Wyspia\'nskiego 27,
         50-370 Wroc{\l}aw, Poland}}}

\date{}\maketitle
\maketitle
\begin{flushleft}\rule[0.4cm]{12cm}{0.3pt}
\parbox[b]{12cm}{\small
\bf Abstract\rm \paragraph{ } A characterization  of an
$h$-hemiregular hemiring  in terms of a fuzzy $h$-ideal  is
provided. Some properties of prime fuzzy $h$-ideals of
$h$-hemiregular hemirings are investigated. It is proved  that a
fuzzy subset $\zeta$ of a hemiring $S$ is a prime fuzzy left (right)
$h$-ideal of $S$ if and only if $\zeta$ is two-valued, $\zeta(0) =
1$, and the set of all $x$ in $S$ such that $\zeta(x) = 1$ is a
prime (left) right $h$-ideal of $S$. Finally, the  similar
properties   for maximal fuzzy left (right) $h$-ideals of
hemirings are considered.\\ \\
{\it Keywords:} Hemiring; $H$-hemiregular hemiring; Fuzzy $h$-ideal;
Prime fuzzy $h$-ideal; Normal fuzzy $h$-ideal \\ \\{\it 2000
Mathematics Subject Classification:} 16Y60; 13E05; 03G25 }
\rule{12cm}{0.3pt}
\end{flushleft}
\footnote{* Corresponding author.\\ \it E-mail address:\rm \ \
zhanjianming@hotmail.com (J. Zhan), dudek@im.pwr.wroc.pl (W.A.
Dudek). }

\subsection*{1. Introduction }

\paragraph{ }There are many concepts of universal algebras generalizing an
associative ring $(R,+,\cdot)$. Some of them -- in particular,
nearrings and several kinds of semirings -- have been proven very
useful. Nearrings arise from rings by cancelling either the axioms
of left or those of right distributivity. The second type of those
algebras $(S,+,\cdot)$, called semirings (and sometimes halfrings),
share the same properties as a ring except that $(S,+)$ is assumed
to be a semigroup rather than a commutative group. Semirings,
ordered semirings and hemirings appear in a natural manner in some
applications to the theory of automata and formal languages (see
\cite{Aho}). It is a well known result that regular languages form
so-called star semirings. According to the well known theorem of
Kleene, the languages, or sets of words, recognized by finite-state
automata are precisely those that are obtained from letters of input
alphabets by the application of the operations sum (union), product,
and star (Kleene closure). If a language is represented as a formal
series with the coefficients in a Boolean semiring, then the Kleene
theorem can be well described by the Sch\"utzenberger Representation
theorem. Moreover, if the coefficient semiring is a field, then  the
corresponding syntactic algebra of the series (see   \cite{Wech} for
details ) has a finite rank if and only if the series are rational.

Many-valued logic   has been proposed to model phenomena in which
uncertainty and vagueness are involved. One of the most general
classes of the many-valued logic   is the BL-logic defined as the
logic of continuous $t$-norms. But in fact, BL-logic  is a
commutative lattice-ordered semiring. So, {\L}ukasiewicz logic,
G\"odel logic and Product logic, as special cases of BL-logic, are
special cases of semirings.

The class of $K$-fuzzy semirings $(K\cup \{+\infty\}, {\rm min},
{\rm max})$, where $K$ denotes a subset of the power set of $R$
which is closed  under the operations ${\rm min}$, $+$, or ${\rm
max}$, has many interesting applications. Min-max-plus computations
(and suitable semirings) are used in several areas, e.g., in
   differential equations. Continuous  timed Petri nets can be
modelled by using generalized polynomial recurrent equations in the
``(min,+) semiring" (see \cite{CGQ}). It is interesting to observe
that the fuzzy calculus, which is used for artificial intelligence
purposes, indeed involves essentially ``(min,max) semirings" (see
\cite{DP} for more details and references). Moreover, the same
hemirings can be used to study  fundamental concepts of the automata
theory such as nondeterminism (cf. \cite{Sim2}). Many other
applications with references can be found in a guide to the
literature on semirings and their applications \cite{Gla}.

Ideals of semirings play a central role in the structure theory and
  are useful for many purposes. However,  they do not in
general coincide with the usual ring ideals if $S$ is a ring and,
for this reason, their use is somewhat limited in trying to obtain
analogues of ring theorems for semirings. Indeed, many results in
rings apparently have no analogues in semirings using only ideals.
M.Henriksen defined in \cite{12} a more restricted class of ideals
in semirings, which is called the class of $k$-ideals, with the
property that if the semiring $S$ is a ring then a complex in $S$ is
a $k$-ideal if and only if it is a ring ideal. Another more
restricted class of ideals has been given in hemirings by K.Iizuka
\cite{13}. However, in an additively commutative semiring $S$,
ideals of a semiring coincide with ``ideals" of a ring, provided
that the semiring is a hemiring. We now call this ideal an $h$-ideal
of the hemiring $S$. The properties of $h$-ideals and also
  $k$-ideals of hemirings were thoroughly investigated by D.R. La
Torre in \cite{22} and by using the $h$-ideals and $k$-ideals, D.R.
La Torre established some analogous ring theorems for hemirings.
Other important results were obtained in \cite{Bo},\cite{Sg},[14-17]
and [22]. Recently, Y.B.Jun \cite{16} considered the fuzzy setting
of $h$-ideals of hemirings.

In this paper, we introduce the concept  of $h$-hemiregularity as a
generalization of the regularity in rings. Next we describe prime
fuzzy $h$-ideals of hemirings and characterize prime fuzzy
$h$-ideals of $h$-hemiregular hemirings by fuzzy $h$-ideals.
Finally, we investigate the properties of normal and maximal fuzzy
left $h$-ideals of hemirings.

\subsection*{2. Preliminaries }

\paragraph{ } Recall that a {\it semiring} is an algebraic system $(S,+,\cdot)$
consisting of a non-empty set $S$ together with two binary
operations on $S$ called addition and multiplication (denoted in the
usual manner) such that $(S,+)$ and $(S,\cdot)$ are semigroups and
the following distributive laws

\vspace{2mm}\centerline{$a(b+c)=ab+ac$ \ \ and \ \ $(a+b)c=ac+bc$}

\vspace{2mm}\noindent are satisfied for all $a,b,c\in S$.

By {\it zero} of a semiring $(S,+,\cdot)$ we mean an element $0\in
S$ such that $0\cdot x=x\cdot 0=0$ and $0+x=x+0=x$ for all $x\in
S$. A semiring with zero and a commutative semigroup $(S,+)$ is
called a {\it hemiring}.

\paragraph{Example 2.1.}(i) A simple example of an infinite hemiring is the set of all
non-negative integers with usual addition and multiplication.

(ii) The set $S=\{0, 1, 2, 3\}$ with the following Cayley tables:
\begin{center}
\begin{tabular}{c|cccc}
          + & 0 & 1 & 2 & 3  \\ \hline
          0 & 0 & 1 & 2 & 3  \\
          1 & 1 & 1 & 2 & 3  \\
          2 & 2 & 2 & 2 & 3  \\
          3 & 3 & 3 & 3 & 2
   \end{tabular}
\ \ \ \ \ \ \ \ \
\begin{tabular}{c|cccc}
          . & 0 & 1 & 2 & 3  \\ \hline
          0 & 0 & 0 & 0 & 0  \\
          1 & 0 & 1 & 1 & 1  \\
          2 & 0 & 1 & 1 & 1  \\
          3 & 0 & 1 & 1 & 1
\end{tabular}
  \end{center}
is   a finite hemiring.

\paragraph{ }A {\it left ideal} of a semiring is a subset $A$ of $S$ closed
with respect to the addition and such that $SA\subseteq A$. A left
ideal $A$ of $S$ is called a {\it left $k$-ideal} if for any $y,z\in
A$ and $x\in S$ from $x+y=z$ it follows $x\in A$. A {\it right
$k$-ideal} is defined analogously.

A left ideal $A$ of a hemiring $S$ is called a {\it left $h$-ideal}
if for any $x,z\in S$ and $a,b\in A$ from $x+a+z=b+z$ it follows
$x\in A$. A {\it right $h$-ideal} is defined analogously. Every left
(respectively, right) $h$-ideal is a left (respectively, right)
$k$-ideal. The converse is not true (cf. \cite{22}).

A fuzzy set $\mu$ of a semiring $S$ is called a {\it fuzzy left
ideal } if for all $x,y\in S$ we have

\begin{enumerate}
\item[$(F_1)$] \ $\mu(x+y)\ge\min\{\mu(x), \mu(y)\}$,
\item[$(F_2)$] \ $\mu(xy)\ge\mu(y)$.
\end{enumerate}

\paragraph{ }Note that a fuzzy left ideal $\mu$ of a semiring $S$ with zero $0$
satisfies also the inequality $\mu(0)\ge\mu(x)$ for all $x\in S$.

\paragraph{Definition 2.2.}A  fuzzy left ideal $\mu$ of $S$ is called a {\it fuzzy left
$k$-ideal } if for all $x,y,z\in S$
\[x+y=z\longrightarrow \mu(x)\ge\min\{\mu(y),\mu(z)\}.\]

A {\it fuzzy right $k$-ideal} is defined analogously. The basic
properties of fuzzy $k$-ideals in semirings are described by Baik
and Kim in \cite{2}.

\paragraph{Definition 2.3.}A fuzzy left ideal $\mu$ of a hemiring $S$ is called a
{\it fuzzy left $h$-ideal } if for all $a,b,x,z\in S$
 \[x+a+z=b+z\longrightarrow \mu(x)\ge\min\{\mu(a),\mu(b)\}.\]

\paragraph{ }A {\it fuzzy right $h$-ideal} is defined similarly. Of course, every
fuzzy left (respectively, right) $h$-ideal is a fuzzy left
(respectively, right) $k$-ideal. The converse is not true (cf.
[12]).

From the Transfer Principle in fuzzy set theory, cf.  \cite{KD}, it
follows that a fuzzy set $\mu$ defined on $X$ can be characterized
by level subsets, i.e. by sets of the form
\[
U(\mu;t)=\{x\in X\,|\,\mu(x)\ge t\},
 \]
where $t\in[0,1]$. Namely, as it is proved in \cite{KD}, for any
algebraic system $\frak{A}=(X,\mathbb{F})$, where $\mathbb{F}$ is a
family of operations (also partial) defined on $X$, the Transfer
Principle can be formulated in the following way:

\paragraph{Lemma 2.4.}{\it A fuzzy set $\mu$ defined on $\frak{A}$ has the property $\mathcal
P$ if and only if all non-empty level subsets $U(\mu;t)$ have the
property $\mathcal{P}$.}

 \paragraph{ }For example, a fuzzy set $\mu$ of a semiring
$S$ is a fuzzy left  ideal if and only if all non-empty level
subsets $U(\mu;t)$ are left ideals of $S$. Similarly, a fuzzy set
$\mu$ in a hemiring $S$ is a fuzzy left $h$-ideal of $S$ if and only
if each non-empty level subset $U(\mu;t)$ is a left $h$-ideal of
$S$.

As a simple consequence of the above property, we obtain the
following proposition, which was first proved in \cite{16}.

\paragraph{Proposition 2.5.}{\it Let $A$ be a
non-empty subset of a hemiring $S$. Then a fuzzy set $\mu_A$ defined
by
\[
\mu_A(x)=\left\{\begin{array}{ll}
 t& {\rm if } \ x\in A\\[2pt]
 s& {\rm otherwise}
 \end{array}\right.
 \]
where $0\le s<t\le 1$, is a fuzzy left $h$-ideal of $S$ if and only
if $A$ is a left $h$-ideal of $S$.}

\paragraph{Definition 2.6.} Let $\mu$ and $\nu$ be fuzzy sets in a hemiring $S$. Then
the $h$-product of $\mu$ and $\nu$ is defined by
\[
(\mu\circ_h\nu)(x)=\sup\limits_{x+a_1b_1+z=a_2b_2+z}
\left(\min\{\mu(a_1),\mu(a_2),\nu(b_1),\nu(b_2)\}\right)
\]
and $(\mu\circ_h\nu)(x)=0$ if $x$ cannot be expressed as
$x+a_1b_1+z=a_2b_2+z$.

\paragraph{Lemma 2.7.} {\it If $\mu$ and  $\nu$ are fuzzy left
$h$-ideals in a hemiring $S$, then so is $\mu\cap\nu$, where
$\mu\cap\nu$ is  defined by
\[(\mu\cap\nu)(x)=\min\{\mu(x),\nu(x)\}
\]
for all $x\in S$. Moreover, if $\mu$ and $\nu$  are  a fuzzy right
$h$-ideal and a fuzzy left $h$-ideal, respectively,  then
$\mu\circ_h \nu\subseteq \mu\cap\nu$.}

\subsection*{3. $H$-hemiregularity }

\paragraph{Definition 3.1.} A hemiring $S$ is said to be {\it $h$-hemiregular} if for each $a\in
S$, there exist $x_1, x_2, z\in S$ such that $a+ax_1a+z=ax_2a+z$.

\paragraph{ }It is not difficult to observe that in the case of rings the
$h$-hemiregularity coincides with the classical regularity of rings.

\paragraph{Example 3.2.} Let $S$ be the set of all non-negative integers $N_0$ with an
element $\infty$ such that $\infty\ge x$ for all $x\in N_0$.
Consider two operations: $a+b=\max\{a,b\}$ and $a\cdot
b=\min\{a,b\}$. It is easy to check that $(S,+,\cdot)$ is an
$h$-hemiregular hemiring.

\paragraph{ }The {\it $h$-closure} $\overline{A}$ of $A$ in a hemiring $S$ is
defined as
\[
\overline{A}=\{x\in S\,|\,x+a_1+z=a_2+z\;\;{\rm for\;some}\;\;
a_1,a_2\in A, \ z\in S\}.
 \]
It is clear that  if $A$ is a left ideal of $S$, then $\overline{A}$
is the smallest left $h$-ideal of $S$ containing $A$. We also have
$\overline{\overline{A}}=\overline{A}$ for each $A\subseteq S$.
Moreover, $A\subseteq B\subseteq S$ implies $\overline{A}\subseteq
\overline{B}$.

\paragraph{Lemma 3.3.} {\it Let $S$ be a hemiring and $A,
B\subseteq S$, then $\overline{AB}=\overline{\overline{A}\
\overline{B}}$.}

\paragraph{Proof.} Because $A\subseteq \overline{A}$ and
$B\subseteq\overline{B}$, then $AB\subseteq \overline{A}\
\overline{B}$, and in   consequence,
$\overline{AB}\subseteq\overline{\overline{A}\ \overline{B}}$.

To prove the converse inclusion, let $x\in\overline{A}$ and
$y\in\overline{B}$. Then there exist $a_i\in A$, $b_i\in B$ and
$z_1,z_2\in S$ such that $x+a_1+z_1=a_2+z_1$ and
$y+b_1+z_2=b_2+z_2$. Putting $z=xz_2+2a_1z_2+z_1y+2z_1b_1+z_1z_2$,
we see that $z\in S$ and

$a_2b_2+a_1b_1+z $

$=a_2b_2+a_1b_1+xz_2+2a_1z_2+z_1y+2z_1b_1+z_1z_2$

$=a_2b_2+xz_2+a_1z_2+z_1y+z_1b_1+z_1z_2+a_1b_1+z_1b_1+a_1z_2$

$=a_2b_2+xz_2+a_1z_2+z_1(y+b_1+z_2)+a_1b_1+z_1b_1+a_1z_2$

$=a_2b_2+xz_2+a_1z_2+z_1(z_2+b_2)+a_1b_1+z_1b_1+a_1z_2$

$=a_2b_2+xz_2+a_1z_2+z_1z_2+z_1b_2+a_1b_1+z_1b_1+a_1z_2$

$=a_2b_2+(x+a_1+z_1)z_2+z_1b_2+a_1b_1+z_1b_1+a_1z_2$

$=a_2b_2+(a_2+z_1)z_2+z_1b_2+a_1b_1+z_1b_1+a_1z_2$

$=a_2b_2+a_2z_2+z_1b_2+z_1z_2+a_1b_1+z_1b_1+a_1z_2$

$=a_2(b_2+z_2)+z_1(b_2+z_2)+a_1b_1+z_1b_1+a_1z_2$

$=(a_2+z_1)(b_2+z_2)+a_1b_1+z_1b_1+a_1z_2$

$=(x+a_1+z_1)(y+b_1+z_2)+a_1b_1+z_1b_1+a_1z_2$

$=xy+a_1b_1+xb_1+a_1y+xz_2+a_1z_2+z_1y+z_1b_1+z_1z_2+a_1b_1+z_1b_1+a_1z_2$

$=xy+(a_1b_1+xb_1+z_1b_1)+(a_1y+a_1b_1+a_1z_2)+(xz_2+a_1z_2+z_1y+z_1b_1+z_1z_2)$

$=xy+(a_1+x+z_1)b_1+a_1(y+b_1+z_2)+(xz_2+a_1z_2+z_1y+z_1b_1+z_1z_2)$

$=xy+(a_2+z_1)b_1+a_1(b_2+z_2)+(xz_2+a_1z_2+z_1y+z_1b_1+z_1z_2)$

$=xy+a_2b_1+a_1b_2+xz_2+2a_1z_2+z_1y+2z_1b_1+z_1z_2$

$=xy+a_2b_1+a_1b_2+z$.

So, $a_2b_2+a_1b_1+z =xy+a_2b_1+a_1b_2+z$, whence we can deduce
$xy\in\overline{AB}$ because $a_ib_j\in AB$ and $z\in S$. This means
that $xy\in\overline{AB}$ for $x\in\overline{A}$,
$y\in\overline{B}$.

Now let $z^{\prime}\in\overline{A}\ \overline{B}$ be arbitrary. Then
$z^{\prime}=\sum\limits_{i=1}^n x_iy_i$ for some $x_i\in
\overline{A}$ and $y_i\in \overline{B}$. Thus
$z^{\prime}\in\overline{AB}$, i.e. $\overline{A}\
\overline{B}\subseteq \overline{AB}$, whence
$\overline{\overline{A}\ \overline{B}}\subseteq
\overline{\overline{AB}}=\overline{AB}$. Therefore
$\overline{\overline{A}\ \overline{B}}= \overline{AB}$. \ \ $\Box$

\paragraph{Lemma 3.4.} {\it If $A$ and $B$ are, respectively,
right and left $h$-ideals of a hemiring $S$, then
$\overline{AB}\subseteq A\cap B$.}

\paragraph{Proof.} Let $x\in \overline{AB}$, then $x+
\sum\limits_{i=1}^m a_ib_i+z= \sum\limits_{j=1}^n a_j' b_j'+z$ for
$a_i,  a_j'\in A$, $b_i, b_j'\in B$ and $z\in S$. Since $A$ is a
right  $h$-ideal of $S$ and $(S,+)$ is a commutative semigroup,
elements $\,\sum\limits_{i=1}^m a_ib_i$, $\sum\limits_{j=1}^n
a_j'b_j'$ are in $A$, and in consequence, $x\in A$. Similarly, we
can prove that $x\in B$. So, $x\in A\cap B$, i.e.
$\overline{AB}\subseteq A\cap B$.\ \ $\Box$

\paragraph{Lemma 3.5.} {\it A hemiring $S$ is $h$-hemiregular if and only if
for any right $h$-ideal $A$ and any left $h$-ideal $B$ we have
$\overline{AB}=A\cap B$.}

\paragraph{Proof.} Assume that $S$ is $h$-hemiregular and $a\in
A\cap B$. Then there exist $x_1,x_2,z\in S$ such that
$a+ax_1a+z=ax_2a+z$. Since $A $ is a right $h$-ideal of $S$, we have
$ax_i\in A$ and $ax_ia\in AB$ for $i=1,2$. Thus $a\in
\overline{AB}$, which implies $A\cap B\subseteq \overline{AB}$.
This, by Lemma 3.4, gives $\overline{AB} =A\cap B$.

Conversely, let $a\in S$. Then, as it is not difficult to verify,
$aS+N_0a$, where $N_0=\{0,1,2,\ldots\}$, is the principal right
ideal of $S$ generated by $a$. Consequently,
$\overline{(aS+N_0a)}$ is a right $h$-ideal of $S$. Therefore
 \[
\overline{(aS+N_0a)}=\overline{(aS+N_0a)}\cap S=
\overline{\overline{(aS+N_0a)}S}= \overline{(aS+N_0a)S}=
\overline{aS}
 \]
because $S$ is trivially an $h$-ideal of itself. Thus
 \[
a=a\cdot 0+1\cdot a\in aS+N_0a\subseteq\overline{(aS+N_0a)}
=\overline{aS}.
 \]
Similarly, $a\in \overline{Sa}$. Hence
 \[
a\in \overline{aS}\cap\overline{Sa}=
\overline{\overline{aS}\cdot\overline{Sa}}=
\overline{aSSa}\subseteq\overline{aSa},
 \]
since $\overline{aS}$ and $\overline{Sa}$ are, respectively, right
and left $h$-ideals of $S$. This shows that there exist $x_1,x_2
,z\in S$ such that $a+ax_1a+z=ax_2a+z$. So, $S$ is a $h$-hemiregular
hemiring.\ \ $\Box$

\paragraph{ }Now, we characterize $h$-hemiregular hemirings by fuzzy
$h$-ideals.

\paragraph{Theorem 3.6.}{\it A hemiring $S$ is $h$-hemiregular if and only if for any fuzzy right
$ h$-ideal $\mu$ and fuzzy left $h$-ideal $\nu$ we have
$\mu\circ_h\nu=\mu\cap\nu$.}

\paragraph{Proof.}Let $S$ be an $h$-hemiregular hemiring. Then
$\,\mu\circ_h\nu\subseteq\mu\cap\nu$ by Lemma 2.7. For any $a\in S$
there exist $x_1,x_2,z\in S$ such that $a+ax_1a+z=ax_2a+z$. Thus
\[\arraycolsep=.5mm
\begin{array}{rl}
(\mu\circ_h\nu)(a)&=\sup\limits_{a+ax_1a+z=ax_2a+z}
\left(\min\{\mu(ax_1),\mu(ax_2),\nu(a)\}\right)\\
&\ge\min\{\mu(ax_1),\mu(ax_2),\nu(a)\}\\[2pt]
&\ge\min\{\mu(a),\nu(a)\}\\[2pt]
&=(\mu\cap\nu)(a),
\end{array}
 \]
i.e. $\,\mu\cap\nu\subseteq \mu\circ_h\nu$, whence
$\,\mu\circ_h\nu=\mu\cap\nu$.

Conversely, let $C$ and $D $ be, respectively right and left
$h$-ideals  of $S$. Then, as it is not difficult to check (cf.
\cite{16}), their characteristic functions $\chi_C$ and $\chi_D$
are, respectively, fuzzy right $h$-ideal and fuzzy left $h$-ideal.
Moreover, by Lemma 3.4, $\overline{CD}\subseteq C\cap D$. Let $a\in
C\cap D$. Then $\chi_C(a)=1=\chi_D(a)$. Thus
 \[
 (\chi_C\circ_h\chi_D)(a)=(\chi_C\cap\chi_D)(a)=\min\{\chi_C(a),\chi_D(a)\}=1.
 \]
So, $\min\{\chi_C(a_1),\chi_D(b_1),\chi_C(a_2),\chi_D(b_2)\}=1$ for
some $a_1,a_2, b_1,b_2$ satisfying the equality
$a+a_1b_1+z=a_2b_2+z$. But then $\chi_C(a_i)=1=\chi_D(b_i)$ for
$i=1,2$, which implies $a_i\in C $ and $b_i\in D$. Therefore
$a\in\overline{CD}$. Hence $\overline{CD}= C\cap D $. Lemma 3.5
completes the proof.\ \ $\Box$

\subsection*{4. Prime fuzzy left $h$-ideals }

\paragraph{ }A prime left (right) $h$-ideals of hemirings are defined in the
same way as prime ideals in rings, i.e. we say that left (right)
$h$-ideal $P$ of a hemiring $S$ is {\it prime} if $P\ne S$ and for
any two left (right) $h$-ideals $A$ and $B$ of $S$ from $AB\subseteq
P$ it follows either $A\subseteq P$ or $B\subseteq P$.

\paragraph{Definition 4.1.} A fuzzy left (right) $h$-ideal $\zeta$ of a hemiring $S$
is said to be {\it prime} if $\zeta$ is a non-constant function and
for any two fuzzy left (right) $h$-ideals $\mu$ and $\nu $ of $S$,
$\mu\circ_h\nu\subseteq\zeta$ implies $\mu\subseteq\zeta$ or
$\nu\subseteq\zeta$.

\paragraph{Example 4.2.} The fuzzy subset
\[
\mu(n)=\left\{\begin{array}{ll}1 & \mbox{ if $n$ is even
 }\\[2pt]
0.2 & \mbox{ otherwise}\end{array}\right.
 \]
defined on the set $N_0$ of all non-negative integers is a prime
fuzzy left $h$-ideal of a hemiring $(N_0,+,\cdot)$.

\paragraph{ }Using the Transfer Principle  (cf. Lemma 2.4) we can
prove the following:

\paragraph{Proposition 4.3.} {\it A fuzzy set $\chi_P$ of a hemiring $S$ is a prime fuzzy
left $($right$)$ $h$-ideal of $S$ if and only if $P$ is a prime left
$($right$)$ $h$-ideal of $S$.}

\paragraph{Theorem 4.4.} {\it A fuzzy subset $\zeta$ of a hemiring $S$ is a
prime fuzzy left $($right$)$ $h$-ideal of $S$ if and only if
\begin{enumerate}
\item[$(i)$] \ $\zeta^0=\{x\in S\,|\,\zeta(x)=\zeta(0)\}$ is a prime
left $($right$)$ $h$-ideal of $S$,
\item[$(ii)$] \ ${\rm Im}\,\zeta=\{\zeta(x)\,|\,x\in S\}$ contains exactly two
elements,
\item[$(iii)$] \ $\zeta(0)=1$.
\end{enumerate}}

\paragraph{Proof.} We prove this theorem only for left $h$-ideals. For right
$h$-ideals the proof is very similar.

Let $\zeta$ be a prime fuzzy left $h$-ideal. Then  it is easy to
check that $\zeta^0$ is a prime left $h$-ideal. Suppose that ${\rm
Im}\,\zeta$ has more than two values. Then there exist two
elements $a,b\in S\backslash\zeta^0$ such that
$\zeta(a)\ne\zeta(b)$. Without loss of generality we can assume
that $\zeta(a)<\zeta(b)$. Since $\zeta$ is a fuzzy left $h$-ideal
and $b\not\in\zeta^0$, it follows that
$\zeta(a)<\zeta(b)<\zeta(0)$. So, there exist $r,t\in [0,1]$ such
that
\[
\zeta(a)<r<\zeta(b)<t<\zeta(0).   \ \ \ \ \ \ \ \ (*)
\]

Let $\nu$ and $\omega$ be fuzzy left $h$-ideals defined by $\nu=r
\chi_{\langle a\rangle}$ and $\omega=t \chi_{\langle b\rangle}$,
where $\chi_{\langle a\rangle}$, $\chi_{\langle b\rangle}$ are
characteristic functions of ideals generated by $a$ and $b$,
respectively. Then, for any $x\in S$, which cannot be expressed in
the form $x+a_1b_1+z=a_2b_2+z$, where $z\in S$, $a_1,a_2\in\langle
a\rangle$ and $b_1,b_2\in\langle b\rangle$, we have
$(\nu\circ_h\omega)(x)=0$. Otherwise,
\[
(\nu\circ_h\omega)(x)=\sup\limits_{x+a_1b_1+z=a_2b_2+z}
\left(\min\{\nu(a_1),\nu(a_2),\omega(b_1),\omega(b_2)\}\right)
=\min\{r,t\}=r.
 \]

Since $\zeta$ is a fuzzy left $h$-ideal, from
$x+a_1b_1+z=a_2b_2+z$ it follows that
\[
\zeta(x)\ge\min\{\zeta(a_1b_1),\zeta(a_2b_2)\}\ge
\min\{\zeta(b_1),\zeta(b_2)\}\ge r.
 \]
So, $(\nu\circ_h\omega)(x)\le\zeta(x)$, whence
$\nu\circ_h\omega\subseteq\zeta$, which implies
$\nu\subseteq\zeta$ or $\omega\subseteq\zeta$ because $\zeta$ is a
fuzzy prime left $h$-ideal. Therefore $\nu(a)=r\le\zeta(a)$ or
$\omega(b)=t\le\zeta(b)$ which contradicts to (*). Consequently,
${\rm Im}\,\mu$ contains exactly two elements.

To prove $(iii)$ suppose that $\zeta$ is a prime fuzzy left
$h$-ideal and $\zeta(0)\ne 1$. Then, according to $(ii)$,\
Im$\,\zeta=\{\alpha_1,\alpha_2\}$, where $0\le\alpha_1<\alpha_2<
1$. Since $\zeta(0)=\zeta(0\cdot x)\ge \zeta(x)$ for all $x\in S$,
we have $\zeta(0)=\alpha_2$. Thus
 \[
\zeta(x)=\left\{\begin{array}{ll}\alpha_2 & \mbox{ if } \ x\in\zeta^0\\
\alpha_1 & \mbox{ otherwise. }\end{array}\right.
 \]

Consider, for fixed $a\in\zeta^0$ and $b\in S\setminus\zeta^0$,
two fuzzy subsets
 \[
\mu(x)=\left\{\begin{array}{ll} t & \mbox{ if } \ x\in\langle a\rangle\\
0 & \mbox{ otherwise}\end{array}\right. \ \ \ \ {\rm and } \ \ \ \
\nu(x)=\left\{\begin{array}{ll} r & \mbox{ if } \ x\in\langle b\rangle\\
0 & \mbox{ otherwise,}\end{array}\right.
 \]
where $0\le\alpha_1< r<\alpha_2< t\le 1$.

It is clear that $\mu$ and $\nu$ are fuzzy left $h$-ideals of $S$.

If $x$ does not satisfy the equality $x+a_1b_1+z=a_2b_2+z$, where
$a_1,a_2\in\langle a\rangle$, $b_1,b_2\in\langle b\rangle$ and
$z\in S$, then $(\mu\circ_h\nu)(x)=0$. Otherwise,
 \[
 (\mu\circ_h\nu)(x)=\sup\limits_{x+a_1b_1+z=a_2b_2+z}
\left(\min\{\mu(a_1),\mu(a_2),\nu(b_1),\nu(b_2)\}\right)
=\min\{t,r\}=r.
 \]

By $(i)$, $\zeta^0$ is a prime left $h$-ideal. If
$a_1,a_2\in\langle a\rangle$, then $a_1,a_2\in\zeta^0$ because
$a\in\zeta^0$ and $\langle a\rangle\subseteq\zeta^0$. This implies
$x\in\zeta^0$. Thus $\zeta(x)=\alpha_2>r=(\mu\circ_h\nu)(x)$.
Therefore, $\mu\circ_h\nu\subseteq\zeta$. But
$\mu(a)=t>\alpha_2=\zeta(a)$ and $\nu(b)=r>\alpha_1=\zeta(b)$,
which gives $\mu\nsubseteq\zeta$ and $\nu\nsubseteq\zeta$. This
contradicts to the assumption that $\zeta$ is a prime fuzzy left
$h$-ideal of $S$. Hence $\zeta(0)=1$.

Conversely, assume that the above conditions are satisfied. Then
$\zeta(0)=1$ and Im$\zeta=\{\alpha,1\}$ for some $0\le\alpha<1$.
Moreover, $\zeta(x+y)\ge\min\{\zeta(x),\zeta(y)\}$ for $x,y\in S$
because $\zeta(x+y)<\min\{\zeta(x),\zeta(y)\}$ implies
$\zeta(x)=\zeta(y)=1$, i.e. $x,y\in\zeta^0$, whence $x+y\in\zeta^0$,
and consequently $\zeta(x+y)=1$, which is impossible. Similarly
$\zeta(xy)\ge\zeta(y)$ since $\zeta(y)=1$ implies $xy\in\zeta^0$,
whence $\zeta(xy)=1$. This means that $\zeta$ is a fuzzy left ideal
of $S$. In fact, $\zeta$ is a fuzzy left $h$-ideal. It is prime.
Indeed, if there exist two fuzzy left $h$-ideals
$\mu\nsubseteq\zeta$ and $\nu\nsubseteq\zeta$ such that
$\mu\circ_h\nu\subseteq\zeta$, then $\mu(x_0)>\zeta(x_0)$ and
$\nu(y_0)>\zeta(y_0)$ for some $x_0,y_0\in S$. It is possible only
in the case when $\zeta(x_0)=\zeta(y_0)=\alpha$, i.e. when
$x_0,y_0\notin\zeta^0$. Since $\zeta^0$ is prime, then there exists
$r\in S$ such that $x_0ry_0\notin\zeta^0$. Otherwise,
$x_0Sy_0\subseteq \zeta^0$, whence $(Sx_0)(Sy_0)\subseteq \zeta^0$.
So, $\overline{(Sx_0)(Sy_0)}\subseteq \overline{\zeta^0}=\zeta^0$,
because $\zeta^0$ is a left $h$-ideal of $S$. Moreover,
$\,\overline{Sx_0}\
\overline{Sy_0}\subseteq\overline{\overline{Sx_0}\
\overline{Sy_0}}=\overline{(Sx_0)(Sy_0)}$ by Lemma 3.3. Thus
$\overline{Sx_0}\ \overline{Sy_0}\subseteq\zeta^0$, and consequently
$\overline{Sx_0}\subseteq\zeta^0$ or
$\overline{Sy_0}\subseteq\zeta^0$. In the first case
$\overline{\langle x_0\rangle\langle x_0\rangle }\subseteq
\overline{Sx_0}\subseteq\zeta^0$, whence $\overline{\langle
x_0\rangle }\subseteq\zeta^0$. So, $x_0\in \langle x_0\rangle
\subseteq\overline{\langle x_0\rangle }\subseteq\zeta^0$. This is a
contradiction. Also the second case yields a contradiction.

Let $a=x_0ry_0$. Then $\zeta(a)=\alpha$. Consequently, by the
assumption
\[
(\mu\circ_h\nu)(a)\le\zeta(a)=\alpha.   \ \ \ \ \ \ \ \ (**)
\]

Obviously $a+x_0ry_0=2x_0ry_0$. Thus
$a+x_0(ry_0)+z=(2x_0)(ry_0)+z$ for any $z\in S$. Therefore for
$a=x_0ry_0$ we have
\[\arraycolsep=.5mm
\begin{array}{rl}
(\mu\circ_h\nu)(a)&=\sup\limits_{a+x_1y_1+z=x_2y_2+z}
\left(\min\{\mu(x_1),\mu(x_2),\nu(y_1),\nu(y_2)\}\right)\\
&\ge\min\{\mu(x_0),\mu(2x_0),\nu(ry_0)\}\ge\min\{\mu(x_0),\nu(ry_0) \}\\[2pt]
&\ge\min\{\mu(x_0),\nu(y_0)\}>\alpha ,
\end{array}
 \]
since $\mu(x_0)>\alpha$ and $\nu(y_0)>\alpha$.

This contradicts (**). Hence for any fuzzy left $h$-ideals $\mu$ and
$\nu$ of $S$, $\mu\circ_h\nu\subseteq\zeta$ implies
$\mu\subseteq\zeta$ or $\nu\subseteq\zeta$. This completes the
proof.\ \ $\Box$

\paragraph{Corollary 4.5.} {\it A fuzzy subset $\zeta$ of a ring $R$ is a prime fuzzy
ideal of $R$ if and only if
\begin{enumerate}
\item[$(i)$] \ $\zeta^0$ is a prime ideal of $R$,
\item[$(ii)$] \ Im$\,\zeta$ contains exactly two elements,
\item[$(iii)$] \ $\zeta(0)=1$.
\end{enumerate}}

\paragraph{Proof.} Since rings are special case of hemirings and in the case of
rings $h$-ideals are ideals, the above result is a simple
consequence of the above theorem.\ \ $\Box$

\paragraph{ }Now, we give an example which shows that in Theorem 4.4 the
condition $(iii)$ cannot be omitted.

\paragraph{Example 4.6.} The set $N_0$ of all non-negative integers is a
hemiring with respect to usual addition and multiplication. Consider
the following fuzzy subsets of $N_0$:
\[
\zeta(n)=\left\{\begin{array}{ll}0.5 & \mbox{ if $n$ is even,}\\[2pt]
0.2 & \mbox{ otherwise}\end{array}\right.
 \]

\[
\mu(n)=\left\{\begin{array}{ll}0.7 & \mbox{ if $n$ is even,}\\[2pt]
0 & \mbox{ otherwise}\end{array}\right.
\]

 \[
\nu(n)=\left\{\begin{array}{ll}0.3 & \mbox{ if $n=3k$ for some
$k\in N_0$}\\[2pt]
0 & \mbox{ otherwise}\end{array}\right.
 \]
Then $\zeta$, $\mu$ and $\nu$ are fuzzy left $h$-ideals such that
$\mu\circ_h\nu\subseteq\zeta$. Moreover, $\zeta^0=\langle
2\rangle=2N_0$ is a prime left $h$-ideal, but
$\mu(2)=0.7>0.5=\zeta(2)$ and $\nu(3)=0.3>0.2=\zeta(3)$. Thus,
$\zeta$ is not a prime fuzzy left $h$-ideal of $N_0$.\ \ $\Box$

\subsection*{5. Normal fuzzy left $h$-ideals }

\paragraph{Definition 5.1.} A fuzzy left $h$-ideal $\mu$ of a hemiring
$S$ is said to be {\it  normal} if there exists $x\in S$ such that
$\mu(x)=1$.

\paragraph{ }Note that if $\mu$ is a normal fuzzy left $h$-ideal of $S$, then
$\mu(0)=1$, and hence $\mu$  is normal if and only if $\mu(0)=1$.
Obviously, any fuzzy left $h$-ideal containing some normal left
$h$-ideal is normal. Indeed, if $\nu\subseteq\mu$ and $\nu$ is
normal, then $1=\nu(0)\le\mu(0)$, and so $\mu(0)=1$.

 Directly from Theorem 4.4 we obtain the following:

\paragraph{Proposition 5.2.} {\it Every prime fuzzy $h$-ideal of a hemiring is
normal.}

\paragraph{ }The converse is not true, see the following:

\paragraph{Example 5.3.}
The fuzzy subset
\[
\mu(x)=\left\{\begin{array}{ll}1 & \mbox{ if $x\in\langle 4
\rangle$,
 }\\[2pt]
0.5 & \mbox{\  if $x\in\langle 2\rangle -  \langle 4\rangle$,}\\[2pt]
0  & \mbox{ otherwise}\end{array}\right.
 \]
defined on a hemiring $(N_0,+,\cdot)$, where $N_0$ is the set of all
non-negative integers, is   a normal fuzzy left $h$-ideal of $N_0$,
which is not prime.

\paragraph{Proposition 5.4.} {\it Given a fuzzy left $h$-ideal $\mu$ of a
hemiring $S$, let $\mu^+$ be a fuzzy set in $S$ defined by
$\mu^+(x)= \mu(x)+1-\mu(0)$ for all $x\in S$. Then  $\mu^+$ is a
normal fuzzy left $h$-ideal of $S$ which contains $\mu$.}

\paragraph{Proof.} For all $x,y\in S$ we have $\mu^+(0)=
\mu(0)+1-\mu(0)=1$ and
\[\arraycolsep=.5mm
\begin{array}{rl}
\mu^+(x+y)&=\mu(x+y)+1-\mu(0)\ge\min\{\mu(x),\mu(y)\}+1-\mu(0)\\[3pt]
&=\min\{\mu(x)+1-\mu(0),\mu(y)+1-\mu(0)\}\\[3pt]
&=\min\{\mu^+(x),\mu^+(y)\},
\end{array}
\] which proves $(F_1)$. Similarly,
\[
\mu^+(xy)=\mu(xy)+1-\mu(0)\ge\mu(y)+1-\mu(0)=\mu^+(y),
\]
verifies $(F_2)$. Hence $\mu^+$ is a fuzzy left ideal of $S$.

Now, let $a,b,x,z\in S$ be such that $x+a+z=b+z$. Then
\[\arraycolsep=.5mm
\begin{array}{rl}
\mu^+(x)&=\mu(x)+1-\mu(0)\ge\min\{\mu(a),\mu(b)\}+1-\mu(0)\\[3pt]
&=\min\{\mu(a)+1-\mu(0),\mu(b)+1-\mu(0)\}\\[3pt]
&=\min\{\mu^+(a),\mu^+(b)\}.
\end{array}
\]
Therefore, $\mu^+$ is a normal fuzzy left $h$-ideal of $S$, and
obviously $\mu\subseteq \mu^+$.\ \ $\Box$

\paragraph{Corollary 5.5.} {\it Let $\mu$ and $\mu^+$ be as in
Proposition 5.4. If there exists $x\in S$ such that $\mu^+(x)=0$,
then $\mu(x)=0$.}

\paragraph{ }Of course, for any left $h$-ideal $A$ of $S$, the characteristic
function $\chi_A$ of $A$ is a normal fuzzy left $h$-ideal of $S$. It
is clear that $\mu$ is a normal fuzzy left $h$-ideal of $S$ if and
only if $\mu^+=\mu$.

\paragraph{Proposition 5.6.} {\it If $\mu$ is a fuzzy left $h$-ideal of
$S$, then $(\mu^+)^+=\mu^+$. Moreover, if $\mu$ is normal, then
$(\mu^+)^+=\mu$.}

\paragraph{Proof.}   Straightforward.\ \ $\Box$

\paragraph{Theorem 5.7.} {\it Let $\mu$ be a fuzzy left $h$-ideal of a hemiring
$S$ and let $f:[0, \mu(0)]\rightarrow[0, 1]$ be an increasing
function. Then a fuzzy set $\mu_f: S\rightarrow[0,1]$ defined by
$\mu_f(x)=f(\mu(x))$ is a fuzzy left $h$-ideal of $S$. In
particular, if $f(\mu(0))=1$, then $\mu_f$ is
 normal; if $f(t)\ge t$ for all $t\in [0, \mu(0)]$, then
$\mu\subseteq\mu_f$.}

\paragraph{Proof.}Indeed, for all $x, y\in S$ we have
\[\arraycolsep=.5mm
\begin{array}{rl}
\mu_f(x+y)&=f(\mu(x+y))\ge f(\min\{\mu(x),\mu(y)\})\\[2pt]
&=\min\{f(\mu(x)), f(\mu(y))\}\\[2pt]
 &=\min\{\mu_f(x),\mu_f(y)\},
\end{array}
\]
which proves $(F_1)$. Similarly,
\[
\mu_f(xy)=f(\mu(xy))\ge f(\mu(y))=\mu_f(y), \] which proves
$(F_2)$. Hence $\mu_f$ is a fuzzy left ideal of $S$.

If $a,b,x,z\in S$ are such that $x+a+z=b+z$, then
 \[\arraycolsep=.5mm
 \begin{array}{rl}
\mu_f(x)&=f(\mu(x))\ge f(\min\{\mu(a),
\mu(b)\})\\[2pt]
&=\min\{f(\mu(a)),f(\mu(b))\}\\[2pt]
&=\min\{\mu_f(a),\mu_f(b)\}.
\end{array}
 \]
Therefore $\mu_f$ is a fuzzy left $h$-ideal of $S$.  If
$f(\mu(0))=1$, then $\mu$ is normal. Assume that $f(t)=f(\mu(x))\ge
\mu(x)$ for any $x\in S$, which proves $\mu\subseteq\mu_f$.\ \
$\Box$

\paragraph{ }Let $\mathcal{N}(S)$ denote the set of all normal fuzzy left
$h$-ideals of $S$. Note that $\mathcal{N}(S)$ is a poset under the
set inclusion.

\paragraph{Theorem 5.8.}{\it Let $\mu\in \mathcal{N}(S)$ be non-constant
such that it is a maximal element of $(\mathcal{N}(S),\subseteq )$.
Then $\mu$ takes only two values $0$ and $1$.}

\paragraph{Proof.}  Since $\mu$ is normal, we have $\mu(0)=1$. Let
$\mu(x)\ne 1$ for some $x\in S$. We claim that $\mu(x)=0$. If not,
then there exists $x_0\in S$ such that $0<\mu(x_0)<1$. Define on $S$
a fuzzy set $\nu$ putting $\nu(x)=(\mu(x)+ \mu(x_0))/2$ for all
$x\in S$. Then, clearly $\nu$ is well-defined and for all $x,y\in S$
we have
 \[\arraycolsep=.5mm
\begin{array}{rl}
\nu(x+y)&=\frac{1}{2}(\mu(x+y)+\mu(x_0))\ge\frac{1}{2}(\min\{\mu(x),
\mu(y)\}+\mu(x_0))\\[2pt]
&=\min\{\frac{1}{2}(\mu(x)+\mu(x_0)), \frac{1}{2}(\mu(y)+ \mu(x_0))\}\\[2pt]
&=\min\{\nu(x),\nu(y)\},
\end{array}
 \]
which proves $(F_1)$. Similarly,
 \[
\nu(xy)=\frac{1}{2}(\mu(xy)+\mu(x_0))\ge\frac{1}{2}
(\mu(y)+\mu(x_0))=\nu(y),
 \]
which proves $(F_2)$. Hence $\nu$ is a fuzzy left ideal of $S$.

Moreover, for $a,b,x,z\in S$ be such that $x+a+z=b+z$ we have
 \[\arraycolsep=.5mm
 \begin{array}{rl}
 \nu(x)&=\frac{1}{2}(\mu(x)+\mu(x_0))\\[2pt]
 &\ge\frac{1}{2}(\min\{\mu(a),\mu(b)\}+\mu(x_0))\\[2pt]
&=\min\{\frac{1}{2}(\mu(a)+\mu(x_0)), \frac{1}{2}(
\mu(b)+\mu(x_0))\}\\[2pt]
&=\min\{\nu(a),\nu(b)\}.
 \end{array}
 \]
Therefore, $\nu$ is a fuzzy left $h$-ideal of $S$.  By Proposition
5.4 \ $\nu^+$ is a maximal fuzzy left $h$-ideal of $S$. Note that
 \[\arraycolsep=.5mm
\begin{array}{rl}
\nu^+(x_0)&= \nu(x_0)+1-\nu(0)\\[2pt]
&=\frac{1}{2}(\mu(x_0)+\mu(x_0))+1-
\frac{1}{2}(\mu(0)+\mu(x_0))\\[2pt]
&=\frac{1}{2}(\mu(x_0)+1)=\nu(x_0)
 \end{array}
  \]
and $\nu^+(x_0)<1=\nu^+(0)$. Hence $\nu^+ $ is non-constant, and
$\mu$ is not a maximal element of $\mathcal{N}(S)$. This is a
contradiction.\ \ $\Box$

\paragraph{Definition 5.9.} A non-constant fuzzy left $h$-ideal $\mu$ of $S$ is
called {\it maximal} if $\mu^+$ is a maximal element of
$\mathcal{N}(S)$.

\paragraph{Theorem 5.10.} {\it If a fuzzy left $h$-ideal $\mu$ of $S$ is maximal, then

$(i)$ \ $\mu$ is normal,

$(ii)$ \ $\mu$ takes only the values $0$ and $1$,

$(iii)$ \ $\chi_{\mu^0}=\mu$,

$(iv)$ \ $\mu^0$ is a maximal left $h$-ideal of $S$.}

\paragraph{Proof.} Let $\mu$ be a maximal fuzzy left $h$-ideal of $S$. Then
$\mu^+$ is a non-constant maximal element of the poset
$(\mathcal{N}(S),\subseteq )$. It follows from Theorem 5.8 that
$\mu^+$ takes only the values $0$ and $1$. Note that $\mu^+(x)=1$ if
and only if $\mu(x)=\mu(0)$, and $\mu^+(x)=0$ if and only if
$\mu(x)=\mu(0)-1$. By Corollary 5.5, we have $\mu(x)=0$, and so
$\mu(0)=1$. Hence $\mu$ is  normal and $\mu^+=\mu$. This proves
$(i)$ and $(ii)$.

$(iii)$ Obvious.

$(iv)$ It is clear that we can  prove   $\mu^0=\{x\in
S\,|\,\mu(0)=1\}$ is a left $h$-ideal. Obviously $\mu^0\ne S$
because $\mu$ takes two values. Let $A$ be a left $h$-ideal
containing $\mu^0$. Then $\mu_{\mu^0}\subseteq\mu_A$, and in
consequence, $\mu=\mu_{\mu^0}\subseteq\mu_A$. Since $\mu$ is normal,
$\mu_A$ also is normal and takes only two values: $0$ and $1$. But,
by the assumption, $\mu$ is maximal, so $\mu=\mu_A$ or $\mu=\omega$,
where $\omega(x)=1$ for all $x\in S$. In the last case $\mu^0=S$,
which is impossible. So, $\mu=\mu_A$, i.e. $\mu_A=\chi_A$. Therefore
$\mu^0=A$.\ \ $\Box$

\paragraph{Definition 5.11.} A normal  fuzzy left $h$-ideal $\mu$ of a
hemiring $S$ is said to be {\it completely  normal} if there exists
$x\in S$ such that $\mu(x)=0$.

\paragraph{ }Denote by $\mathcal{C}(S)$ the set of all completely normal  fuzzy
left $h$-ideals of $S$. We note that $\mathcal{C}(S)
\subseteq\mathcal{N}(S)$ and the restriction of the partial ordering
$\subseteq$ of $\mathcal{N}(S)$ gives a partial ordering of
$\mathcal{C}(S)$.

\paragraph{Proposition 5.12.} {\it Any non-constant maximal element of
$(\mathcal{N}(S),\subseteq)$ is also a maximal element of
$(\mathcal{C}(S),\subseteq)$.}

\paragraph{Proof.} Let $\mu$ be a non-constant maximal element of
$(\mathcal{N}(S),\subseteq)$. By Theorem 5.10, $\mu$ takes only the
values 0 and 1, and so $\mu(0)=1$ and $\mu(x)=0$ for some $x\in S$.
Hence $\mu\in \mathcal{C}(S)$. Assume that there exists
$\nu\in\mathcal{C}(S)$ such that $\mu\subseteq\nu$. It follows that
$\mu\subseteq\nu$ in $\mathcal{N}(S)$. Since $\mu$ is maximal in
$(\mathcal{N}(S),\subseteq)$ and $\nu$ is non-constant, therefore
$\mu=\nu$. Thus $\mu$ is maximal element of
$(\mathcal{C}(S),\subseteq)$, which ends the proof.\ \ $\Box$

\paragraph{Theorem 5.13.} {\it Every maximal fuzzy left $h$-ideal of a
hemiring $S$ is completely normal.}

\paragraph{Proof.} Let $\mu$ be a maximal  fuzzy left
$h$-ideal of $S$. Then by Theorem 5.10, $\mu$ is normal and
$\mu=\mu^+$ takes only the values 0 and 1. Since $\mu$ is
non-constant, it follows that $\mu(0)=1$ and $\mu(x)=0$ for some
$x\in S$. Hence $\mu$ is completely  normal, which ends the proof.\
\ $\Box$

\subsection*{6. Conclusions }

\paragraph{ } In the present paper, we showed that the basic results of fuzzy sets in
hemirings are similar, but not identical, with the corresponding
results for semirings. So, it is important for us to  study
different types of fuzzy sets  among hemirings, semirings and rings.
In our opinion the future study of fuzzy sets in hemirings and
semirings can be connected with (1) investigating
  semiprime fuzzy $h$-ideals; (2) establishing a fuzzy spectrum of
a hemiring; (3) finding intuitionistic and/or interval-valued fuzzy
sets and triangular norms.  The obtained results can be   used to
solve some social networks problems and to decide whether the
corresponding graph is balanced or clusterable.

\subsection*{Acknowledgements }

\paragraph{ }The authors are highly grateful to referees and Professor Witold
Pedrycz, Editor-in-Chief, for their valuable comments and
suggestions for improving the paper. The authors also express their
thanks to  Professor John Mordeson and language editor for improving
English language.

The research of the first author is partially supported by   the
National Natural Science Foundation of China (60474022)   and  the
Key Science Foundation of Education Committee of Hubei Province,
China ( 2004Z002; D200529001).

\end{document}